\documentclass[12pt]{amsart}

\usepackage[utf8]{inputenc}
\usepackage[T1]{fontenc}
\usepackage{newtxtext}
\usepackage[a4paper, margin=1in]{geometry}
\usepackage{setspace}
\usepackage{amsmath, amssymb, amsthm, amsfonts, mathtools}
\usepackage{graphicx}
\usepackage{float}
\usepackage{microtype}
\usepackage{pgfplots}
\pgfplotsset{width=10cm,compat=1.18}
\usepackage{subcaption}
\usepackage{paralist}
\usepackage{amsrefs}  

\theoremstyle{plain}
\newtheorem{theorem}{Theorem}[section]
\newtheorem{proposition}[theorem]{Proposition}
\newtheorem{lemma}[theorem]{Lemma}

\theoremstyle{definition}
\newtheorem{definition}[theorem]{Definition}
\newtheorem{example}[theorem]{Example}

\theoremstyle{remark}

\DeclarePairedDelimiter\abs{\lvert}{\rvert}
\DeclarePairedDelimiter\norm{\lVert}{\rVert}

\begin{document}

\title[Piecewise Convex  Maps with Countably many Branches]{Quasi-compactness of Frobenius-Perron Operator for Piecewise Convex  Maps with Countably many Branches}

\author{Pawe\l \ G\'ora}
\author{Aparna Rajput}
\address{Department of Mathematiocs and Statistics, Concordia University, Montreal, H3G 1M8, QC, Canada}
\email{pawel.gora@concordia.ca}
\email{a\_ajpu@live.concordia.ca}


\thanks{The research for this publication was supported by NSERC (National Science and Engineering Research Council of Canada) grant,
number RGPIN-2020-06788}
\subjclass[2020]{Primary 37A05; Secondary 37E05}

\date{}

\dedicatory{}

\commby{Pawe\l \ G\'ora}

\begin{abstract}
In this paper, we prove the quasi-compactness of the Frobenius-Perron operator for a piecewise convex map with a countably infinite number of branches on the interval $I=[0,1]$. We establish that for sufficiently high iterates of the map, the iterations are piecewise expanding. By adapting the Lasota-Yorke Inequality to meet the assumptions of the Ionescu-Tulcea and Marinescu ergodic theorem, we demonstrate the existence of an absolutely continuous invariant measure (ACIM) for the map, the exactness of the dynamical system and the quasi-compactness of the Frobenius-Perron operator induced by the map. These findings imply a multitude of strong ergodic properties for the system.
\end{abstract}

\maketitle

\section{\textbf{Introduction}}
\label{sec1}
This paper investigates the existence of absolutely continuous invariant measures (ACIMs) for a class $\mathcal{T}$ of piecewise convex maps with a countably infinite number of branches, defined on the unit interval $[0,1]$. Understanding ACIMs is crucial for analyzing the long-term behavior and chaotic nature of deterministic dynamical systems.

To formalize this study, we consider the measure space  \((I,\mathcal{B},m )\), where $I=[0,1], \mathcal{B}$ is the Borel $\sigma$ algebra of subsets of $I$, and $m$ is the Lebesgue measure. By $L^1$ we denote the space of integrable functions with respect to the Lebesgue measure $m$. All transformations, $\tau: I \rightarrow I$, we consider, are piecewise strictly monotonic and hence non-singular transformation (i.e. $
m(\tau^{-1}(A)) = 0 \quad \text{whenever} \quad m(A) = 0, \quad \forall A \in \mathcal{B}
$). A measure $\mu$ on $\mathcal{B}$ is $\tau$-invariant if it remains unchanged under the action of $\tau$, i.e., $\mu(\tau^{-1} A)=\mu(A)$,  for all $A \in \mathcal{B}$. To study ACIMs for piecewise monotonic, non-singular transformations, we use the Frobenius-Perron operator \( P_\tau \) induced by \(\tau\), which acts on  \( L^1 \) space and, for piecewise monotonic maps, can be represented as
\begin{equation*}
    P_\tau f(x) =\sum_{y\in \tau^{-1}\{x\}} \frac{f(y)}{\left|\tau'(y)\right|}.
\end{equation*}
The Frobenius-Perron operator \( P_\tau \) has key properties such as linearity, positivity, contractivity, and preservation of integrals. A measure of the form \( f \cdot m \), where \( f \) is a density function in \( L^1 \) and \( m \) is the Lebesgue measure on \(I\), is \(\tau\)-invariant if and only if \( P_\tau f = f \). This means \( f \) is a fixed point of \( P_\tau \). These properties make \( P_\tau \) essential for studying ACIMs, as it allows us to characterize invariant measures in terms of fixed points of the operator $P_\tau$. For more information about ACIMs, the Frobenius-Perron operator, and their mutual relationship, we refer the reader to  \cite{boyarsky1997} or \cite{lasota1994}.

The piecewise convex maps considered until recently had a finite number of branches. We recall the standard assumptions.
Let $I=[0,1]$, we say that $\tau$ belongs to the class of piecewise convex maps $ \mathcal{T}^{pc}$ if it satisfies the following conditions :
\begin{enumerate}[(1)]
\item There exist a finite partition $0=a_0<a_1<.....<a_n=1$ such that $\tau|_{[a_{i-1},a_i)}$ is continuous, strictly increasing and convex for $i=1,2...n$.\\
\item $\tau(a_{i-1})=0$ and $\tau'(a_{i-1})>0$ for $i=1,2,...n$.\\
\item $\tau'(0)=1/ \alpha>1$.
\end{enumerate}
 The first to study such maps were Lasota and Yorke \cites{lasota1973existence}. They discovered the three important properties, restated in Theorem \ref{Thm 2.1.6} below, and used them to prove the existence of the ACIM $\mu$ and the exactness of the system $(\tau,\mu)$. Additional properties of piecewise convex maps were shown in \cites{jablonski1976,jablonski1983,jablonski1983b,jablonski1984,jablonski1991}. Generalizations, weakenings of the assumptions, and random maps were studied in \cites{bose2002,inoue1992,inoue2024}.
Recently, in \cites{rahman2019, gora_submit, islam_submit} Lasota \& Yorke's results were generalized to the case when a piecewise convex map has a countably infinite number of branches. The Lasota-Yorke type \mbox{inequality} they proved is for the pair $(L^\infty, L^1)$ and only for non-increasing functions. We prove the \mbox{Lasota-Yorke} type inequality for the pair $(BV_I, L^1)$, which allows us to show the quasi-compactness of the operator $P_\tau$ induced by $\tau$.
This fully describes the behavior of the system $(\tau,\mu)$ and implies several strong
ergodic properties, in particular, the exponential decay of correlation. These results could not be obtained by employing the previously used methods.

In Section \ref{sec2}, we define piecewise convex maps with a countably infinite number of branches. These maps are defined on the partition of $I=[0,1]$ into disjoint open subintervals ${I_i = (a_i, b_i)}_{i=1}^\infty$, whose complement has Lebesgue measure zero. Each restriction $\tau_i$ to $I_i$ is an increasing, convex, differentiable function, with $\sum_{i \ge 1} \frac{1}{\tau_i'(a_i)} < +\infty$ and $\tau'(0) > 1$, if $0$ is not a limit point of partition endpoints. We denote the class of such maps by $\mathcal{T}$. We study the Frobenius-Perron operator $P_\tau$ associated with these maps. We prove several key properties of $P_\tau$, including its effect on non-increasing functions and bounds on its norm. 
We show that if $\tau$ belongs to the class $\mathcal{T}$, its iterates $\tau^n$ retain the piecewise convex structure and summability condition on derivatives. We demonstrate that the set of preimages of partition points is dense in $[0,1]$, and the derivatives of sufficiently high iterates are uniformly bounded below by a constant greater than one, indicating piecewise expanding behavior.

 In Section \ref{section: 2.2},  we recall basic facts about piecewise expanding maps with a countably infinite number of branches. These maps are defined on a partition of $I$ into disjoint subintervals ${I_i = (a_i, b_i)}_{i=1}^\infty$, each homeomorphically mapped by $\tau$ onto its image. The expansion behavior is quantified by $g(x)=\frac{1}{\abs{\tau'(x)}}$, with $g(x) \leq \beta < 1$, defining the class $\mathcal{T}_E$ . We prove that if $\tau$ is in $\mathcal{T}$, then its iterates $\tau^n$ are in $\mathcal{T}_E$ for all large enough $n$.  We prove Theorem \ref{thm 2.2.9}, which guarantees the existence and uniqueness of a normalized absolutely continuous invariant measure for $\tau$ and shows that the system is exact. We study the Frobenius-Perron operator $P_\tau$ and its action on functions in $BV_I$, where $BV_I = \{ f \in L^1 : \norm{f}_{BV} < +\infty \}$ and follow  \cites{keller1985,rychlik1983} to establish a Lasota-Yorke inequality for the functions in $BV_I$ and apply Ionescu-Tulcea and Marinescu theorem  to prove the quasi-compactness of $P_\tau$ and its implications for ACIMs (Theorem \ref{Theorem 2.2.7}).
 
 Random number generators are popular applications of piecewise monotonic maps. If the distribution we want to simulate has a decreasing density, then an appropriate piecewise
convex map would be a natural choice. For example if the density is $g=e^{1-x}/(e-1)$ any of the maps $h^{-1}\circ\tau_k\circ h$, where $h=(e/(e-1))(1-e^{-x})$ and $\tau_k(x)=kx $ (mod 1), $k=2,3,...$, would work.
\begin{figure}[H]
  \centering
  \begin{subfigure}[b]{0.55\textwidth}
    \includegraphics[width=\textwidth]{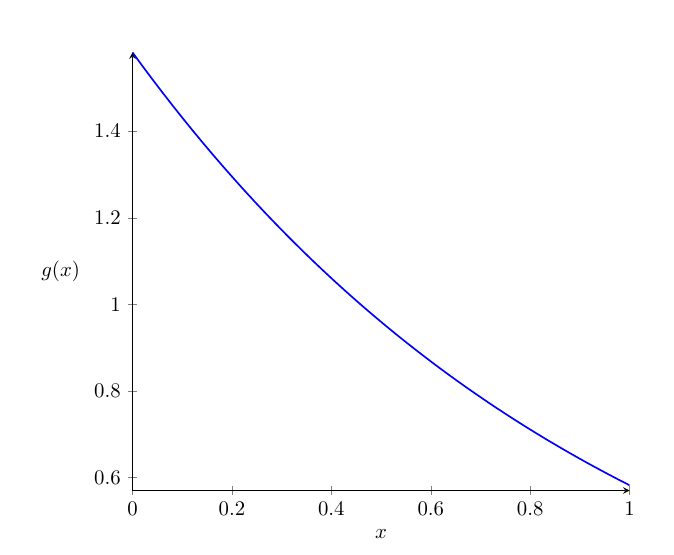}
    \caption{Plot of density function $g(x)=\frac{e^{1-x}}{(e-1)}$.}
    \label{Intro image1}
  \end{subfigure}
  \hfill
  \begin{subfigure}[b]{0.42\textwidth}
    \includegraphics[width=\textwidth]{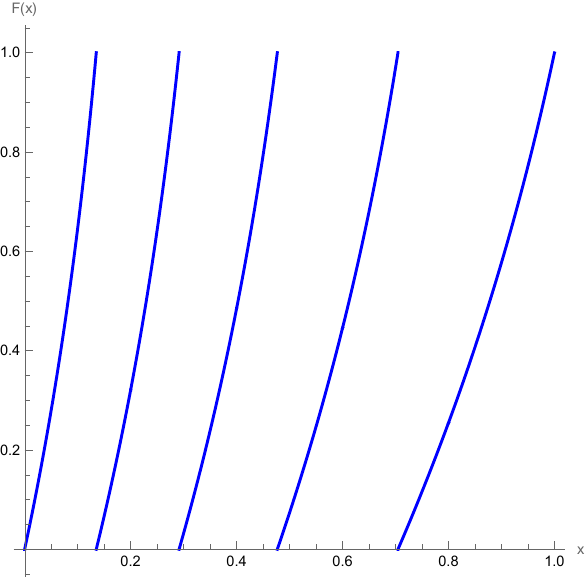}
     \caption{Plot of $F(x)=h^{-1}\circ\tau_k\circ h$, $k=5$, a piecewise convex map preserving the measure with density $g$.}
    \label{Intro image2}
  \end{subfigure}
\end{figure}
In some dynamical systems, a map with infinitely many branches can appear when we consider the first return map of a map with only a few branches. For example, suppose \( \tau \) is a piecewise convex map on \( [0,1] \) where the derivative at zero satisfies \( \tau'(0) = 1 \). This lack of expansion near zero makes it impossible to use standard techniques to show that \( \tau \) has an ACIM. However, if we define the first return map \( \tau_{\rm fr} \), which shows how points return to a specific part of the interval that avoids zero (e.g., \( [\epsilon, 1] \)), we obtain a new map. This map \( \tau_{\rm fr} \) is piecewise convex with infinitely many branches. Crucially, each branch exhibits expansion at its left endpoint, enabling us to apply our results to establish the existence of an ACIM for the  map \( \tau_{\rm fr}\) and, if additional conditions are satisfied, for the original map $\tau$.
\section{ Piecewise convex map with countably infinite number of branches}
\label{sec2}
Let $I=[0,1]$, consider a measure space $(I, \mathcal{B}, m)$, where $m$ is the Lebesgue measure on $I$ and $\mathcal{B}$ is the Borel $\sigma$-algebra on $I$. We define the class of piecewise convex map with countably infinite number of branches :
\begin{definition} \label{def 2.1.1} Let $I=[0,1]$ and let $\mathcal P=\{I_i=(a_i,b_i)\}_{i=1}^\infty$ be a countably infinite family of open disjoint subintervals of $I$ such that Lebesgue measure of $I\setminus\bigcup_{i\ge 1}I_i$ is zero.
We define a piecewise convex map $\tau:I\rightarrow I$ on the partition $\mathcal P$ as follows:\\
(2.1) \label{2.1} For $i=1,2,3..., \tau_i=\tau_{|I_i}$ is an increasing, convex, and differentiable function with $\lim_{x\to a_i^+}\tau_i(x)=0$.
We define $\tau_i(a_i)=0$ and $\tau_i(b_i)=\lim_{x\to b_i^-}\tau_i(x)$. The values $\tau_i'(a_i)$ are also defined by (right-) continuity.\\
(2.2) \label{2.2} We assume $$\sum_{i\ge 1}\frac 1{\tau_i'(a_i)} < +\infty.$$\\
(2.3) \label{2.3}  If $x=0$ is not a limit point of the partition points, then we have $\tau'(0)=1/\alpha>1$, for some $0<\alpha<1$.\\
 We will denote the set of  maps satisfying conditions (2.1)-(2.3) by $\mathcal T$.
\end{definition}
\begin{example}\label{Exam 2.1.2} 
   Let us consider a  map $\tau:I\rightarrow I$ with countable number of branches defined on $I_i=(a_i , a_{i+1})$, with $a_i=\frac{i-1}{i}$ for $i=1,2, \dots$ by,
\begin{equation*}
    \tau_i(x)=\left(\frac{x-a_i}{a_{i+1} -a_i}\right) .
\end{equation*}
Since $\tau_i$ is linear, it is a convex function on each subinterval $I_i$. It is also differentiable on each $I_i$, with the derivative given by 
 \begin{equation*}
     \tau_i'(x)=\frac{1}{a_{i+1} -a_i} .
 \end{equation*}
 Note that the function $\tau_i$ takes the value $0$ at the left end of the interval $I_i$, for each $i=1,2,...$
  \begin{figure}[H]
     \centering
     \includegraphics[width=0.85\linewidth]{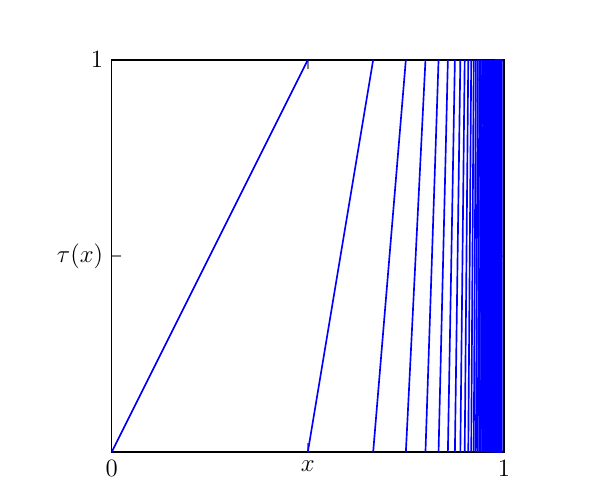}
     \caption{A piecewise convex map of example \ref{Exam 2.1.2}.}
 \end{figure} 
 \end{example}
 To check the summability condition (2.2), note that 
 \begin{equation*}
      \tau_i'(a_i)=\frac{1}{a_{i+1} -a_i}=i(i+1) .
 \end{equation*}
 Hence
 \begin{equation*}
     \sum_{i\ge 1}\frac 1{\tau_i'(a_i)}=\sum_{i\ge 1}\frac {1}{i(i+1)} \leq \sum_{i\ge 1} \frac{1}{i^2},
 \end{equation*}
  which is a convergent series. Therefore, the condition $(2.2)$ is satisfied.\\
  Notice that $a_1=0$ is not a limit point of the partition end points, $\tau_1'(a_1)=2$.
 Hence, this example of a piecewise convex map satisfies all the conditions of Definition \ref{def 2.1.1}.
\begin{lemma}\label{lem 2.1.3} Let $\tau \in \mathcal T$, and assume that $x = 0$ is a limit point of $\{a_i\}$. Then we have $lim_{x\to 0^+}\tau'(x) = +\infty$. Moreover, for some $0 < r < 1$ and some $\alpha<1$, we have
$$\sum_{a_i<r}\frac 1{\tau'(a_i)} =\alpha <1 .$$
\end{lemma}
\begin{proof}
Condition (2.2) implies that for any $M>0$ the inequality $\tau'(a_i)\le M$  can be satisfied only for a finite number
of points $a_i$. This implies the first claim of the lemma. The second claim follows by the fact that the tails of a convergent series
converge to $0$.
\end{proof}
We now consider an example similar to Example \ref{Exam 2.1.2}, but with $I_i=\bigl( \frac{1}{i+1}, \frac{1}{i}\bigl)$; note that $0$ is a limit point of $\{b_i\}$.
\begin{example}\label{Exam 2.1.4}
    Let us consider a map $\tau:I \rightarrow I$ with countable number of branches defined on $ I_i=(b_{i+1},b_i)$, with $b_i=\frac{1}{i}$ for $i=1,2, \dots$ by 
\begin{equation*}
     \tau_i(x)=\left(\frac{x-b_{i+1}}{b_{i} -b_{i+1}}\right) . 
\end{equation*}
\end{example}

\begin{figure}[H]
    \centering
    \includegraphics[width=.75\linewidth]{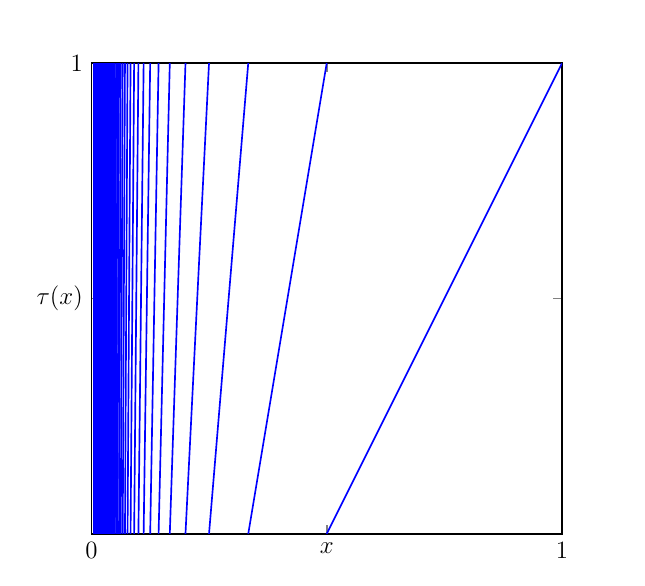}
    \caption{$\tau$ with a countable number of branches, with $x=0$ a limit point of $\{b_i\}$.}
    \label{fig:2.2}
\end{figure}
 \bigskip
 The Frobenius-Perron operator induced on $L^1_m$ by the map $\tau\in\mathcal T$ is
\begin{equation}
\begin{aligned}
P_\tau f (x)=\sum_{i\ge 1} \frac {f(\tau_i^{-1}(x))}{\tau_i'(\tau_i^{-1}(x))}\chi_{|\tau_i(I_i)} (x).
\end{aligned}
\end{equation}
The properties proven in the following theorem have been shown and used before to prove the existence of ACIMs for \( \tau \in \mathcal{T}^{\text{pc}} \). Here we again establish these properties again, but for \( \tau \in \mathcal{T} \). Because of the differences between \( \mathcal{T} \) and \( \mathcal{T}^{\text{pc}} \), this requires slightly different methods.
The following theorem summarizes the properties of the maps in $\mathcal T$.
\begin{theorem} \label{Thm 2.1.6} Let $\tau \in \mathcal T$, and let $f:I\to \mathbb R^+$ be a non-increasing function. Then:

(1) $P_\tau f$ is a non-increasing function;\\

(2) for any $x\in I \setminus \{0\}$ we have $f(x)\le \frac 1x \| f\|_1 $;\\

(3) $\| P_\tau f\|_\infty\le \alpha \| f\|_\infty + D \| f\|_1$, where $0<\alpha<1$ is the number specified in Definition \ref{def 2.1.1} and Lemma \ref{lem 2.1.3}, and $D>0$ is a constant.\\
\end{theorem}
\begin{proof}
 (1) Since $\tau_i$ is increasing, $\tau_i^{-1}$ is also increasing on $\tau_i(I_i)$. We prove that $P_\tau f$ is non-increasing. Let $f\in L^1$ be non-increasing. Then for any $x\leq y$ we have $f(x)\geq f(y)$. We fix $i \geq 1$. We need to show that
  \begin{equation*}
      \frac{f(\tau_i^{-1}(x))}{\tau_i'(\tau_i^{-1}(x))} \cdot \chi_{\tau_i (I_i)} (x)\geq \frac{f(\tau_i^{-1}(y))}{\tau_i'(\tau_i^{-1}(y))} \cdot \chi_{\tau_i(I_i)} (y).
  \end{equation*}
  For a non increasing function $f\in L^1$ and $x,y\in \tau_i(I_i)$, we have,
  \begin{equation*}
  \begin{aligned}
     \tau_i^{-1}(x) < \tau_i^{-1}(y) 
      \implies f( \tau_i^{-1}(x))>f( \tau_i^{-1}(y)).
  \end{aligned}
  \end{equation*}
  By convexity of $\tau_i$ we know that $\tau_i'$ is increasing, which gives
  \begin{equation*}
      \begin{aligned}
          \tau_i'( \tau_i^{-1}(x))<\tau_i'( \tau_i^{-1}(y))
          \implies \frac{1}{\tau_i'( \tau_i^{-1}(x))}>\frac{1}{\tau_i'( \tau_i^{-1}(y))}.
      \end{aligned}
  \end{equation*}
 Recall that the characteristic function $\chi_A (x)$ is $1$ if $x\in A$ and is $0$ if $x\notin A$. If $x\in \tau_i(I_i)$ and $y\notin \tau_i(I_i) $ then $\chi_{\tau_i(I_i)} (x)\geq\chi_{\tau_i(I_i)} (y)$. The case \( x \notin \tau_i(I_i) \), \( y \in \tau_i(I_i) \) is impossible because \( \tau_i(I_i) = [0, \tau_i(b_i)] \), so if \( y \in \tau_i(I_i) \) and \( x \leq y \), then \( x \in \tau_i(I_i) \) as well. If \( x \notin \tau_i(I_i) \) and \( y \notin \tau_i(I_i) \) or if \( x \in \tau_i(I_i) \) and \( y \in \tau_i(I_i) \) then $\chi_{\tau_i(I_i)} (x)=\chi_{\tau_i(I_i)} (y)$. So for any $x,y\in I$, $x\le y$, we have
  \begin{equation*}
      \chi_{\tau_i(I_i)} (x)\geq\chi_{\tau_i(I_i)} (y).
  \end{equation*}
  Since the product of non-negative non-decreasing functions is non-decreasing, we get that for any $i$,
  \begin{equation*}
      \frac{f(\tau_i^{-1}(x))}{\tau_i'(\tau_i^{-1}(x))} \cdot \chi_{\tau_i(I_i)} (x)\geq \frac{f(\tau_i^{-1}(y))}{\tau_i'(\tau_i^{-1}(y))} \cdot \chi_{\tau_i(I_i)} (y).
  \end{equation*}
  This shows, that $P_\tau f(x)$ is a non-increasing function.\\
  (2)  For every \( x \in I \), since \( f > 0 \) and \( f \) is non-increasing,
\[
\norm{f}_1 = \int_I f \, dm \geq \int_0^x f(t) \, dt \geq x \cdot f(x).
\]
\begin{figure}[H]
      \centering
      \includegraphics[width=0.75\linewidth]{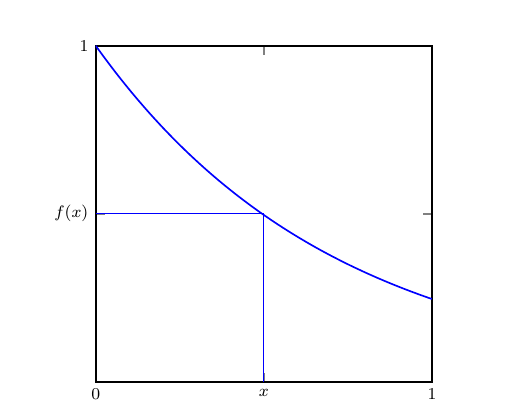}
      \caption{Inequality for (2).}
  \end{figure}
  
(3) \textbf{Case I}: Let $x=0$ be a limit point of $\{a_i\}$. Since $P_\tau f$ is non-increasing, $\|P_\tau f \|_\infty \leq P_\tau f(0)$, and we have, with $r, \alpha$ as in Lemma \ref{lem 2.1.3}, 
   \begin{align}
       P_\tau f(0)&=\sum_{i\geq 1}\frac{f(\tau_i^{-1}(0))}{\tau_i'(\tau_i^{-1}(0))} =\sum_{i\geq 1}\frac{f(\tau_i^{-1}(\tau_i(a_i))}{\tau_i'(\tau_i^{-1}(\tau(a_i))}=\sum_{i\geq 1} \frac{f(a_i)}{\tau_i'(a_i)}\nonumber&\\
        &= \sum_{i: a_i<r} \frac{f(a_i)}{\tau_i'(a_i)}+\sum_{i: a_i\geq r} \frac{f(a_i)}{\tau_i'(a_i)}\leq  \alpha \cdot \|f \|_\infty+\sum_{i: a_i\geq r} \frac{f(a_i)}{\tau_i'(a_i)} \nonumber&\\
        &\leq \alpha \cdot \|f \|_\infty+\sum_{i: a_i\geq r} \frac{1}{a_i \cdot \tau_i'(a_i)}\cdot \| f\|_1.
    \end{align}
For this case, we define $\displaystyle D=\sum_{i: a_i>r} \frac{1}{a_i \cdot \tau_i'(a_i)}
   < \infty$.\\
   \textbf{Case II}: Let  $x=0$ not be a limit point of  $\{a_i\}$. Then $a_i=0$ for some $i$, assume w.l.o.g. that $a_1=0$. Again, since $P_\tau f$ is non-increasing, $\|P_\tau f\|_\infty \leq P_\tau f(0)$, and we have
       \begin{align}
           P_\tau f(0)&=\sum_{i\geq 1}\frac{f(\tau_i^{-1}(0))}{\tau_i'(\tau_i^{-1}(0))}=\sum_{i\geq 1}\frac{f(\tau_i^{-1}(\tau_i(a_i))}{\tau_i'(\tau_i^{-1}(\tau(a_i))}=\sum_{i\geq 1} \frac{f(a_i)}{\tau_i'(a_i)}\nonumber&\\
           &=\frac{f(a_1)}{\tau_1'(a_1)}+\sum_{i\geq 2} \frac{f(a_i)}{\tau_i'(a_i)}\leq \frac{f(0)}{\tau_1'(0)}+\sum_{i\geq 2} \frac{1}{a_i \cdot \tau_i'(a_i)}\cdot \|f\|_1\nonumber&\\
           &\leq \alpha \cdot \|f\|_\infty + D \cdot \|f\|_1.
       \end{align}
For this case, we define $\displaystyle D =\sum_{i\geq 2} \frac{1}{a_i \cdot \tau_i'(a_i)}$, which is finite by condition (2.2).
\end{proof}
Let $\mathcal P^{(n)}=\mathcal P \bigvee \tau^{-1}(\mathcal P)\bigvee\dots\bigvee \tau^{n-1}(\mathcal P) $. We denote the branches of $\tau^n$ by $\tau^{(n)}_i$. Then, $\mathcal P^{(n)}=\left\{I_i^{(n)}=\left(a_i^{(n)},b_i^{(n)}\right)\right\}_{i=1}^\infty$ is a countably infinite family of open disjoint subintervals of $I$ corresponding to $\tau^n$. We have the following results:
 \begin{proposition} \label{Pro 2.1.7}  Let $\mathcal{P}$ be a partition for $\tau$, and let $\mathcal{P}^{(n)}$ denote the partition for $\tau^n$. If $\tau \in \mathcal{T}$ then $\tau^n \in \mathcal{T}$ as well, i.e.,
 \begin{enumerate}
  \item[ \textbf{(a)}]  $\tau^n$  is piecewise increasing on $\mathcal{P}^{(n)}$;
  \item[\textbf{(b)}] $\tau^n$  is piecewise convex on $\mathcal{P}^{(n)}$;
   \item[\textbf{(c)}]  $\tau^n$  is piecewise differentiable on $\mathcal{P}^{(n)}$;
   
    \item[\textbf{(d)}]$\lim_{x\rightarrow \left(a_i^{(n)}\right)^+} \tau_i^{(n)}(x)=0$ for every $i$; 
    
   \item[\textbf{(e)}] the condition $(2.2)$ holds for $\tau^n$, i.e.,
    \begin{equation*}
    \sum_{i\ge 1}\frac {1}{\left(\tau^{(n)}_i\right)'\left(a_i^{(n)}\right)} < +\infty;
    \end{equation*}
    \item[\textbf{(f)}] if $x=0$ is not a limit point of $\{a
    _i\}$ and if condition $(2.3)$ holds for $\tau$  then it holds for $\tau^n$.
    \end{enumerate}
    \end{proposition}
\begin{proof} The proofs of (a), (b) and (c) are simple, since $\tau$ is piecewise increasing, convex and differentiable on $\mathcal{P}$ and we know that the composition of increasing, convex and differentiable functions is also increasing, convex, and differentiable. Hence (a), (b) and (c) hold for $\tau^n$ on $\mathcal{P}^{(n)}$.\\
    To prove (d), we will use induction. From (2.1) we have $\lim_{x\to a_i^+}\tau_i(x)=0$. We consider for $n=2$ one branch of $\tau^2$. The branch $\tau_k^{(2)}=\tau_j\circ \tau_i$ is defined on $I_{i,j}=I_i\cap\tau_i^{-1}(I_j)=\tau_i^{-1}(I_j)$. Since the left endpoint of $\tau_i(I_i)$ is $0$, if the interval $\tau_i^{-1}(I_j)$ is not empty, then it contains $\tau_i^{-1}(a_j)=a_k^{(2)}$, the left endpoint of $I_k^{(2)}$.\\ 
 We have
 \begin{equation*}
     \begin{aligned}
          \lim_{x\to \left(a_k^{(2)}\right)^+} (\tau_j \circ \tau_i)(x)=\tau_j\left(\tau_i\left(a_k^{(2)}\right)\right)=\tau_j(\tau_i(\tau_i^{-1}(a_j)))=\tau_j(a_j)=0.
     \end{aligned}
 \end{equation*}
 Now, we use induction. We assume that the result holds for $\tau^n$. The map $\tau^n$ has infinitely many branches and each branch satisfies the property (d), i.e.,
 \begin{equation*}
         \lim_{x\to \left(a_i^{(n)}\right)^+} \tau_i^{(n)}(x)=0 .
 \end{equation*}
 For $n+1$, if we consider a $k^{th}$ branch of $\tau^{n+1}$, $\tau_k^{(n+1)}=\tau_j\circ \tau_i^{(n)} $. We have,
 \begin{equation*}
     a_k^{(n+1)}=\left(\tau_i^{(n)}\right)^{-1}(a_j).
 \end{equation*}
 and
 \begin{equation*}
    \lim_{x\rightarrow  \left(a_k^{{(n+1)}}\right)^+ }\tau_j\left(\tau_i^{(n)} \left(x \right)\right)=\tau_j\left(\tau_i^{(n)}\left(a_k^{(n+1)}\right)\right)=\tau_j\left(\tau_i^{(n)}\left(\left(\tau_i^{(n)}\right)^{-1}\left(a_j \right) \right) \right)=\tau_j(a_j)=0.
 \end{equation*}
 \\
(e) Let us assume $\displaystyle \sum_{j\ge 1}\frac 1{\tau_j'\left(a_j\right)}=K$. We will prove that
\begin{equation*}
    \displaystyle \sum_{k\ge 1} \frac 1 {\left(\tau_k^{(n+1)}\right)'\left(a_k^{(n+1)}\right)} \le K^{n+1}.
\end{equation*}
We consider for $n=2$ one branch of $\tau^2$. The branch $\tau_k^{(2)}=\tau_j\circ \tau_i$ is defined on $I_{i,j}=I_i\cap\tau_i^{-1}(I_j)=\tau_i^{-1}(I_j)$. Then,
 \begin{equation*}
 \begin{aligned}
    \sum_{j\geq 1} \sum_{i\ge 1}\frac {1}{\left(\tau_j \circ \tau_i\right)'\left(a_k^{(2)}\right)}&=\sum_{j\geq1} \sum_{i\ge 1}\frac{1}{\tau_j'\left(\tau_i\left(a_k^{(2)}\right)\right) \tau_i'\left(a_k^{(2)}\right)}&\\
    & =\sum_{j\ge1} \sum_{i\ge 1}\frac{1}{\tau_j'(a_j) \tau_i'\left(a_k^{(2)}\right)}&\\
     &\leq \sum_{j\ge1} \frac{1}{\tau_j'(a_j)}\sum_{i\ge 1}\frac{1}{\tau_i'\left(a_k^{(2)}\right)}&\\
     &= K\cdot \sum_{i\ge 1}\frac{1}{\tau_i'\left(a_k^{(2)}\right)}.
     \end{aligned}
 \end{equation*}
 Since $\tau_i$ is convex, which implies that $\tau_i'$ is non-decreasing,  we have $\tau_i'(a_i)\leq \tau_i'\left(a_k^{(2)}\right)$, and 
 \begin{equation*}
     \begin{aligned}
          \sum_{k\ge 1} \frac{1}{\left(\tau_k^{(2)}\right)'\left(a_k^{(2)}\right)}=\sum_{j\geq 1} \sum_{i\ge 1}\frac {1}{\left(\tau_j \circ \tau_i\right)'\left(a_k^{(2)}\right)}&\leq K \cdot \sum_{i\geq 1} \frac{1}{\tau_i'(a_i)}&=K \cdot K <\infty.
     \end{aligned}
 \end{equation*}
 For $n+1$, we consider a $k^{th}$ branch of $\tau^{n+1}$, $\tau_k^{(n+1)}=\tau_j\circ \tau_i^{(n)} $. We have,
 \begin{equation*}
     a_k^{(n+1)}=\left(\tau_i^{(n)}\right)^{-1}(a_j).
 \end{equation*}
 and
 \begin{equation*}
      \begin{aligned}
     \sum_{j\geq 1} \sum_{i\ge 1}\frac {1}{\left(\tau_j \circ \tau_i^{(n)}\right)'\left(a_k^{(n+1)}\right)}&=\sum_{j\geq1} \sum_{i\ge 1}\frac{1}{\tau_j'\left(\tau_i^{(n)}\left(a_k^{(n+1)}\right)\right) \left(\tau_i^{(n)}\right)'\left(a_k^{(n+1)}\right)}&\\
     &=\sum_{j\ge1} \sum_{i\ge 1}\frac{1}{\tau_j'(a_j) \left(\tau_i^{(n)}\right)'\left(a_k^{(n+1)}\right)}&\\
     &\leq \sum_{j\ge1} \frac{1}{\tau_j'(a_j)}\sum_{i\ge 1}\frac{1}{\left(\tau_i^{(n)}\right)'\left(a_k^{(n+1)}\right)}= K\cdot \sum_{i\ge 1}\frac{1}{\left(\tau_i^{(n)}\right)'\left(a_k^{(n+1)}\right)}.
 \end{aligned}
 \end{equation*}
 Since $\left(\tau_i^{(n)}\right)'$ is non-decreasing, so $\left(\tau_i^{(n)}\right)'\left(a_i^{(n)}\right)\leq \left(\tau_i^{(n)}\right)'\left(a_k^{(n+1)}\right)$, and we have
 \begin{equation*}
     \begin{aligned}
          \sum_{j\geq 1} \sum_{i\ge 1}\frac {1}{\left(\tau_j \circ \tau_i^{(n)}\right)'\left(a_k^{(n+1)}\right)}&=K \cdot \sum_{i\ge1} \frac{1}{\left(\tau_i^{(n)}\right)'\left(a_i^{(n)}\right)}=K\cdot K^n=K^{n+1} <\infty.
     \end{aligned}
 \end{equation*}
By induction the result holds for any n.\\
 (f)  Let  $x=0$ be not a limit point of  $\{a_i\}$. Then, $x=0$ is not a limit point of $\{a_k^{(n)}\}$ either, and  $\tau'(0)=\frac{1}{\alpha}>1$ gives $(\tau^n)'(0)=\frac{1}{\alpha^n}>1$.
 \end{proof}
 \upshape
Let us consider $\tau \in \mathcal{T}$. We now show that the set of preimages of the partition points under all iterations of $\tau$ is dense in $I$. 
 \begin{proposition}\label{Prop 2.1.8}
     Let $\tau \in \mathcal{T}$. Then the set
     \begin{equation*}
         S=\bigcup_{n=0}^\infty \tau^{-n}(\{a_i,b_i : i=1,2, \dots\})
     \end{equation*} 
     is dense in $[0,1]$.
     \end{proposition}
     \begin{proof}
If $x=0$ is not a limit point of the partition endpoints, we assume w.l.o.g that $a_1 =0$. Let
\begin{equation*}
     S=\bigcup_{n=0}^\infty \tau^{-n}(\{a_i,b_i : i=1,2, \dots\}).
     \end{equation*}
We want to prove that $S$ is dense in $[0,1]$. Let us suppose that it was not true. Then, there exists an interval $[x_0,y_0]\subset [0,1]$ such that
\begin{equation*}
    \tau^n([x_0,y_0])\cap\{a_i,b_i : i=1,2, \dots\}=\emptyset \hspace{0.4cm} \text{for all $n=0,1,2,3...$}.
\end{equation*}
Therefore for each $n$, the points $x_n=\tau^n(x_0)$ and $y_n=\tau^n(y_0)$ belong to the same interval $(a_i,b_i)$. Let $x_n,y_n \in (a_k,b_k)$ and $x_n<y_n$ , $k=1,2,3...$. The map $\tau_k$ is defined on $(a_k,b_k)$, $k=1,2,3...$. Let  $\theta_1$ be the angle between the positive half of the $x$-axis and the segment connecting $(0,0)$ to $(x_n,\tau_k(x_n))$. Similarly, let  $\theta_2$ be the angle between the positive half of the $x$-axis and the segment connecting $(0,0)$ to $(y_n,\tau_k(y_n))$. See Figure \ref{Chap2:fig3} for $k=1$ and Figure \ref{Chap2:fig4} for $k>1$.  First, let $k=1$. For $\tau_1=\tau\vert_{(a_1, b_1)}$ we have
\begin{equation*}
    \tan \theta_1 = \frac{\tau_1(x_n)}{x_n} \hspace{0.2cm} \text{and} \hspace{0.2cm} \tan \theta_2 = \frac{\tau_1(y_n)}{y_n}.
\end{equation*}
Further, since $\tau_1$ is increasing, we have
\begin{equation*}
    \tan \theta_2 \geq \tan \theta_1 \implies \frac{\tau_1(y_n)}{y_n} \geq  \frac{\tau_1(x_n)}{x_n}
\end{equation*}
\begin{equation*}
    \implies \frac{y_{n+1}}{x_{n+1}} =\frac{\tau^{n+1}(y_0)}{\tau^{n+1}(x_0)}=\frac{\tau_1(\tau^n(y_0)}{\tau_1(\tau^n(x_0)} \geq \frac{y_n}{x_n}.
\end{equation*}
\begin{figure}[H]
    \centering
    \includegraphics[width=1\linewidth]{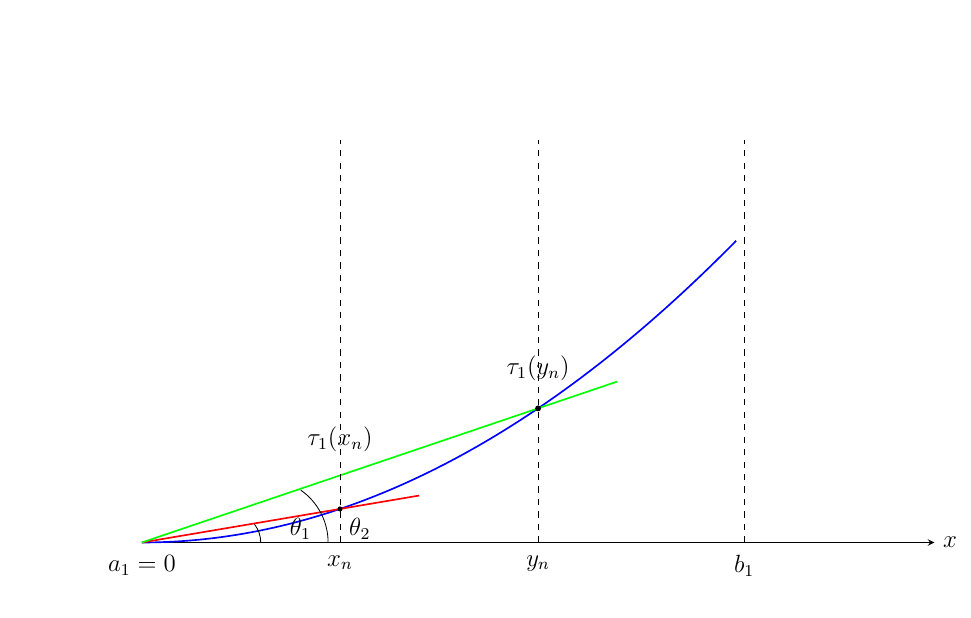}
    \caption{When $0$ is not a limit point of $\{a_i\}$, and $k=1$.}
    \label{Chap2:fig3}
\end{figure}
\begin{figure}[H]
    \centering
    \includegraphics[width=1\linewidth]{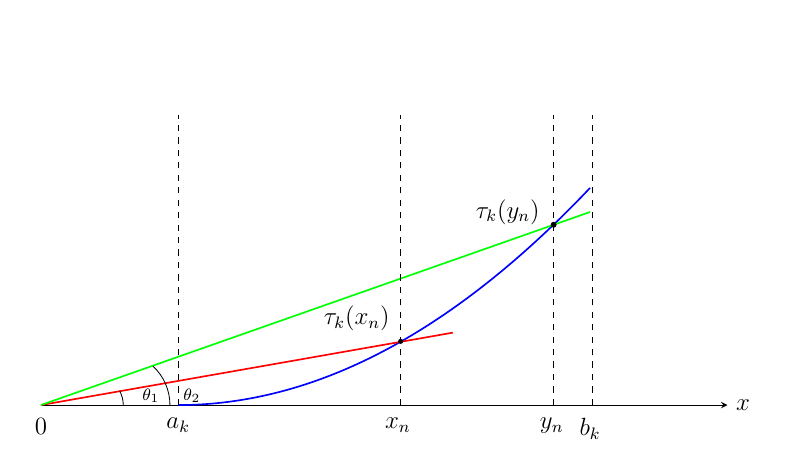}
    \caption{When $0$ is not a limit point of $\{a_i\}$, and $k>1$.}
    \label{Chap2:fig4}
\end{figure}
Now, for $k>1$, $\tau_k=\tau\vert_{(a_k, b_k)}$ we again have
\begin{equation*}
    \tan\theta_1=\frac{\tau_k(x_n)}{x_n} \hspace{0.5cm} \text{and} \hspace{0.5cm} \tan\theta_2=\frac{\tau_k(y_n)}{y_n}.
\end{equation*}
Since $\tau_k$ is increasing on $(a_k,b_k)$ we have,
    \begin{equation}
        \tan\theta_2\geq \tan\theta_1  \implies \frac{\tau_k(y_n)}{y_n} \geq \frac{\tau_k(x_n)}{x_n} \implies \frac{\tau_k(y_n)}{\tau_k(x_n)} \geq \frac{y_n}{x_n}.
    \end{equation}
or,
    \begin{equation*}
        \frac{x_{n+1}}{y_{n+1}}=\frac{\tau_k(x_n)}{\tau_k(y_n)}\leq \frac{x_n}{y_n}.
    \end{equation*}
 Since this holds for $k=1,2,3...$ we obtain for all $n\geq 1$,
    \begin{equation}\label{(2.1.5)}
        \frac{x_{n+1}}{y_{n+1}} \leq \frac{x_n}{y_n}\leq...\leq \frac{x_0}{y_0}.
        \end{equation}
        If $0$ is the limit of the partition points we set $b_E=b_1$. Otherwise, we set $b_E=b_j$ such that $b_j<r$, where $r>0$ is as in Lemma \ref{lem 2.1.3}. We define $I_E=[0,b_E]$. The map $\tau$ is piecewise expanding on  $I_E$.  

Now, we will prove a stronger estimate, valid only for pairs $\{x_n,y_n\}$ such that  $x_n,y_n\in (a_k,b_k)$ with $b_E<a_k$. 
        \begin{figure}[H]
            \centering
            \includegraphics[width=1\linewidth]{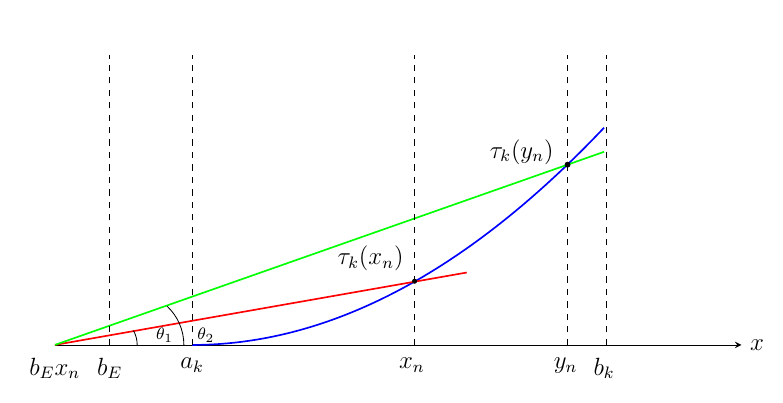}
            \caption{When $0$ is not a limit point of $\{a_i\}$ and $k>1$.}
            \label{Chap2:fig5}
        \end{figure}
        We redefine $\theta_1$ as the angle between the positive half of the $x$-axis and the segment connecting $(b_E\cdot x_n, 0)$ to $(x_n,\tau_k(x_n))$. Similarly, let  $\theta_2$ be the angle between the positive half of the $x$-axis and the segment connecting $(b_E\cdot x_n, 0)$ to $(y_n,\tau_k(y_n))$. See Figure
\ref{Chap2:fig5}. We have
\begin{equation*}
    \tan\theta_1=\frac{\tau_k(x_n)}{x_n-b_Ex_n} \hspace{0.5cm} \text{and} \hspace{0.5cm} \tan\theta_2=\frac{\tau_k(y_n)}{y_n-b_Ex_n},
\end{equation*}
which is valid only if $b_E<a_k$. Since $\tau_k$ is increasing  on $(a_k,b_k)$ we have
    \begin{equation*}
        \tan\theta_2\geq \tan\theta_1\implies \frac{\tau_k(y_n)}{y_n-b_Ex_n} \geq \frac{\tau_k(x_n)}{x_n-b_Ex_n}, 
    \end{equation*}
    or
    \begin{equation*}
        \frac{\tau_k(y_n)}{\tau_k(x_n)} \geq \frac{y_n-b_Ex_n}{x_n-b_Ex_n}= \frac{y_n\left(1-b_E\frac{x_n}{y_n}\right)}{x_n(1-b_E)}.
    \end{equation*}
    By Equation (\ref{(2.1.5)}) we obtain,
    \begin{equation*}
        \frac{1-b_E\frac{x_n}{y_n}}{1-b_E}\geq \frac{1-b_E\frac{x_0}{y_0} }{1-b_E}.
    \end{equation*}
  Thus, for $x_n,y_n\in (a_k,b_k) $ with $b_E<a_k$, we obtain 
\begin{align}
    \frac{y_{n+1}}{x_{n+1}}\geq q\frac{y_n}{x_n},\label{(2.1.6)}
\end{align}
where $\displaystyle q=\frac{1-b_Ex_0/y_0}{1-b_E}>1$.\\
Whenever the interval $(x_n,y_n)$ falls inside the interval $I_E$ it is stretched by $\tau$ as long as it stays in $I_E$.  This expansion ensures that the length of \( (x_n, y_n) \) grows with each iteration. Consequently, the images of the interval $(x_n,y_n)$ leave the interval
$I_E$ after a finite number of steps.  Equation (\ref{(2.1.6)}) implies that $\lim_{n\rightarrow \infty} \frac{y_n}{x_n}=\infty$. Since $\lim\sup_n x_n \geq b_E $, we have $\lim\sup_n y_n =\infty$, which is impossible as it contradicts the fact that $y_n \in [0,1]$ for all $n\in \mathbb{N}$.\\
Hence $S$ is dense in $[0,1]$.
\end{proof}
\begin{proposition}\label{Prop 2.1.9}
Let $\tau \in \mathcal{T}$  and let $\mathcal{P}^{(n)} = \{I_i^{(n)} = (a_i^{(n)}, b_i^{(n)})\}_{i=1}^{\infty}$ be the partition of  $[0,1]$ corresponding to $\tau^n$. Then, given any $\epsilon >0$, there exists $\Bar{n} \in \mathbb{N}$ such that for all $n > \Bar{n}$,
\[
\max_i \left( b_i^{(n)} - a_i^{(n)} \right) < \epsilon.
\]
\end{proposition}
\begin{proof}
We suppose that this was not true, i.e., there exists an $\epsilon_0>0$ such that for any $n\ge 1$
we  can find $I_{i_n}^{(n)}=(a_{i_n}^{(n)},b_{i_n}^{(n)})$ with 
$$b_{i_n}^{(n)}-a_{i_n}^{(n)}\geq \epsilon_0.$$
Let $c_{i_n}^{(n)}$ be the mid point of  the interval $I_{i_n}^{(n)}$. Then,
\[
\left( c_{i_n}^{(n)} - \frac{\varepsilon_0}{2}, c_{i_n}^{(n)} + \frac{\varepsilon_0}{2} \right)
\bigcap \left\{a^{(k)}_j, b^{(k)}_j : j \in \mathbb{N}\hspace{0.1cm}, k \leq n \right\}
= \emptyset .
\] 
We used the fact that each partition  interval corresponding to $\tau^n$ lies inside a partition interval corresponding to $\tau^k$ i.e., the partition $\mathcal{P}^{(n)}$ is a refinement of the partition $\mathcal{P}^{(k)}$, for $k\le n$. \\
Let $c$ be a limit point of $c_{i_n}^{(n)}$. Then,
\[
\left( c - \frac{\varepsilon_0}{3}, c + \frac{\varepsilon_0}{3} \right)
\bigcap \left\{a^{(n)}_j, b^{(n)}_j: j, n \in \mathbb{N}\right\}
=
\left( c - \frac{\varepsilon_0}{3}, c + \frac{\varepsilon_0}{3} \right)
\cap S = \emptyset ,
\] 
which contradicts the denseness of the set $S$ proved in Proposition \ref{Prop 2.1.8}.
\end{proof}
\begin{theorem}\label{thm 2.1.10} Let $\tau \in \mathcal{T}$. Then, there exists $n_0 \in \mathbf{N}$ such that $\inf(\tau^n)'\geq 2$ for all $n \geq n_0$.
\end{theorem}
\begin{proof} Recall that $\mathcal{P}^{(n)}=\left\{\left(a_i^{(n)},b_i^{(n)}\right)\right\}_{i=1}^\infty $ is a partition corresponding to $\tau^n$, and the branch of $\tau^n$ defined on the interval $\left(a_i^{(n)},b_i^{(n)}\right)$ is $\tau^{(n)}_i$. We also know from Proposition \ref{Pro 2.1.7} that $\tau^n$ satisfies conditions (2.1)-(2.3) with respect to the partition $\mathcal{P}^{(n)}$.
Consider the set
 \begin{equation*}
 S=\bigcup_{n=0}^\infty \tau^{-n}(\{a_i,b_i : i=1,2, \dots\}).
 \end{equation*}
 In Proposition \ref{Prop 2.1.8} we proved that $S$ is dense in $[0,1]$. Note that $\inf\tau' >0$.
 Similarly as for the proof of Proposition \ref{Prop 2.1.8},  if $x=0$ is not a limit point of $\{a_i\}$, we assume w.l.o.g. that $a_1=0$.
We also define $I_E=[0,b_E]$, where $b_E=b_1$ if $0$ is not a limit point of the partition points and otherwise  $b_E=b_j$ such that $b_j<r$, where $r>0$ is as in Lemma \ref{lem 2.1.3}.
By Proposition \ref{Pro 2.1.7}, with $\epsilon= \frac{b_E}{2} \inf \tau' >0$, there exist $\Bar{n}\in \mathbf{N}$ such that for any $n>\Bar{n}$,
\begin{equation}
    \max_i\left(b_i^{(n)}-a_i^{(n)}\right)< \epsilon.\label{2.1.9}
\end{equation}
This condition ensures that the length of the longest interval in the partition of $I$ for the $n$-th iterate of $\tau$ is less than $\epsilon$.\\
Let $n>\Bar{n}$. Define
\begin{equation*}
    A_n=\tau^{-n}(0,b_E)=\bigcup_{i=1}^\infty (\tau_i^n)^{-1}((0,b_E)),
\end{equation*}
and $B_n=[0,1]\setminus A_n$. We have
\begin{equation}
    \tau^n(x)<b_E  \hspace{0.5cm} \text{if} \hspace{0.5cm} x\in A_n,
\end{equation}
and
\begin{equation}
    \tau^n(x)\geq b_E  \hspace{0.5cm} \text{if} \hspace{0.5cm} x\in B_n .
\end{equation}
\begin{figure}[H]
    \centering
    \includegraphics[width=1\linewidth]{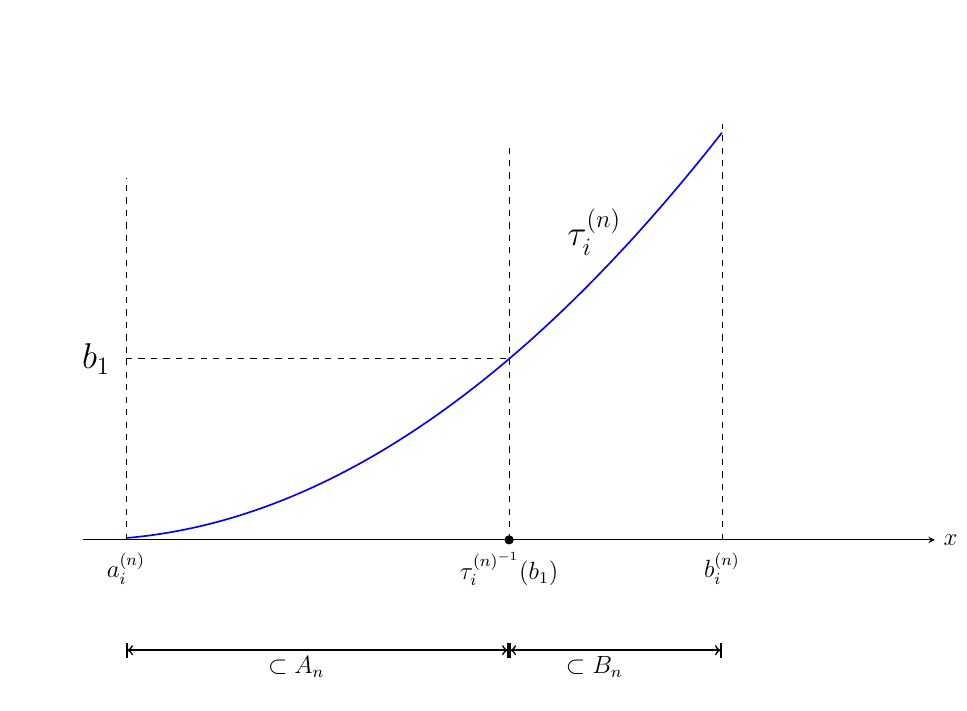}
    \caption{$n$-th iterate of $\tau$ on $I_i^{(n)}$.}
    \label{Chap2:fig6}
\end{figure}
The map $\tau^n$ is increasing on each interval $\left(a_i^{(n)},b_i^{(n)}\right)$, and $\left(\tau_i^{{(n)}}\right)'$ represents the rate of change of $\tau^n$ on this interval. By equation (\ref{2.1.9}), we know that the length of this interval is less than $\epsilon$. Then, the length of the interval $\left(a^{(n)}_i,\left(\tau^{(n)}_i\right)^{-1}(b_E)\right)$ is also less than $\epsilon$. Therefore, $\left(\tau_i^{(n)}\right)'$ must be sufficiently large to ensure that $\tau^n$ increases by at least $b_E$ over an interval of length less than $\epsilon$ which gives,
\begin{equation*}
    \left(\tau_i^{(n)}\right)'(x)\geq \frac{b_E}{\epsilon},
\end{equation*}
for some $x\in\left(a_i^{(n)},\left(\tau^{(n)}_i\right)^{-1}(b_E)\right)=A_n \cap \left(a_i^{(n)},b_i^{(n)}\right)$. Since $\left(\tau_i^{(n)}\right)'$ is increasing, we have the same inequality for all $x\in B_n \cap \left(a^{(n)}_i,b^{(n)}_i\right)$ and hence
\begin{equation}
\begin{aligned}
    \left(\tau \circ \tau_i^{(n)}\right)'(x)&=\tau'\left(\tau_i^{(n)}(x)\right)\cdot \left(\tau_i^{(n)}\right)'\left(x\right) ,&\\
    &\geq \frac{b_E}{\epsilon} \inf \tau'= 2,\label{2.1.12}
    \end{aligned}
\end{equation}
whenever $x\in B_n$ and, $i=1,2,3....$ . For $x\in A_n$ we have
\begin{equation}
\begin{aligned}
   \left(\tau \circ \tau_i^{(n)}\right)'(x)&=\tau'\left(\tau_i^{(n)}(x)\right)\cdot \left(\tau_i^{(n)}\right)'(x)\geq \tau'(0) \left(\tau_i^{(n)}\right)'\left(a_i^{(n)}\right)&\\
    &\geq \tau'(0) \inf(\tau^n)'. \label{2.1.13}
    \end{aligned}
\end{equation}
Inequalities (\ref{2.1.12}) and (\ref{2.1.13}) together give us
\begin{equation*}
    \inf\left(\tau^{n+1}\right)'\geq \min \left({2, \tau'(0) \inf(\tau^n)'}\right),
\end{equation*}
and consequently, by induction we have 
\begin{equation*}
     \inf(\tau^{n})'\geq \min ({2, \tau'(0)^{n-\Bar{n}} \inf(\tau^{\Bar{n}}})'),
\end{equation*}
for all $n>\Bar{n}$. This implies that for sufficiently large $n$ we have $\inf(\tau^n)'\ge 2$.
\end{proof}
\section{ Piecewise expanding map with countably many branches}\label{section: 2.2}
\begin{definition} \label{def 2.2.1} Let $I=[0,1]$, and let $\mathcal P=\{I_i=(a_i,b_i)\}_{i=1}^\infty$ be a countably infinite family of open disjoint subintervals of $I$ such that Lebesgue measure of $I\setminus\bigcup_{i\ge 1}I_i$ is zero. Let $\tau$ be a differentiable map from $\displaystyle \cup_{i\geq 1} I_i$ to the interval $I$, such that for each $i\geq 1$, $\tau_{|I_i}$ extends to a homeomorphism $\tau_i$ of $[a_i,b_i]$ onto its image.\\ Let
\[
g(x) = 
\begin{cases}
    \frac{1}{|\tau_i'(x)|}, & \text{for } x \in I_i, i = 1, 2, \dots ,\\
    0, & \text{elsewhere}.
\end{cases}
\]
We assume $\sup_{x \in I} | g(x) | \leq \beta < 1$, \textrm{var}$(g)< + \infty$. Then, we say $\tau$ is a piecewise expanding map with countably many branches and denote the class of all such $\tau$ by $\mathcal{T}_E$.
\end{definition}
We now prove a lemma about the relation between the families $\mathcal{T}$ and $\mathcal{T}_E$ of maps $\tau: I \rightarrow I$ introduced in Definitions \ref{def 2.1.1} and \ref{def 2.2.1}, respectively.
\begin{theorem}\label{thm 2.2.2} If $\tau\in \mathcal{T}$, then some iterate  $\tau^n\in \mathcal T_E$ for all sufficiently large $n\in \mathbb{N}$.
\end{theorem}
\begin{proof}
   This is a direct consequence of Theorem \ref{thm 2.1.10} and the condition (2.2) in Definition \ref{def 2.1.1}.
    \end{proof}
\normalfont{ Let $\tau\in \mathcal{T}_E$, $\norm{g}_\infty<1$, \textrm{var}$(g)< + \infty$. A piecewise expanding map $\tau$ is non-singular, and the Frobenious-Perron operator corresponding to $\tau$ is
\begin{equation}
    P_\tau f(x)=\sum_{i=1}^\infty \frac{f\left(\tau_i^{-1}(x)\right)}{|\tau'\left(\tau_i^{-1}(x)\right)|} \chi_{\tau(I_i)}(x)=\sum_{y\in\tau^{-1}\{x\}} f(y)g(y).
\end{equation}
\begin{theorem}\label{Prop-2.2.3} Assume \( \tau \in \mathcal{T}_E \). Then, for every sufficiently large \( n \in \mathbb{N} \), there exists \( B_n \in \mathbb{R}^+ \) such that
\[
\textnormal{var}_I P_{\tau^n} f \leq \frac{1}{2} \textnormal{var}_I f + B_n \|f\|_1
\]
for all \( f \in BV_I \).
\end{theorem}
\begin{proof}
    This is a direct consequence of Corollary~3 in Rychlik's paper \cite{rychlik1983}.
\end{proof}

\begin{proposition}\label{prop 2.2.5} 
For every $\tau \in \mathcal{T}$ there exist $n\in \mathbb{N}$ and $C>0$ so that for all $f\in BV_I$, we have
\begin{equation*}
    \|P_{\tau^n} f\|_{BV} \le \frac{1}{2} \cdot \|f\|_{BV} + C \cdot \|f\|_1 .
\end{equation*}
\end{proposition}
\begin{proof}
Recall that for every $f\in BV_I$,
\begin{equation*}
    \|f\|_{BV}= \int_I |f| dm + \textrm{var}_If=\|f\|_1+\textrm{var}_If ,
\end{equation*}
and hence, for every $n\in \mathbb{N}$,
\begin{equation}
\begin{aligned}
    \|P_{\tau^n} f\|_{BV}=\|P_{\tau^n}f\|_1+\textrm{var}_I P_{\tau^n}f
    \leq \|f\|_1+ \textrm{var}_I P_{\tau^n}f.\label{2.2.7}
    \end{aligned}
\end{equation}
With Theorem \ref{thm 2.2.2}, pick $n_0 >1$ so large that $\tau^{n_0}\in \mathcal{T}_E$.
Every $f\in BV_I$ has a version $f^*$ with a minimal variation, characterized by the property that, for all $x_0\in (0,1)$,
\begin{equation*}
    f^*(x_0)\in \left[\lim_{x\rightarrow x_0^-}\ f^*,\lim_{x\rightarrow x_0^+}\ f^*\right].
\end{equation*}
All one-sided limits always exist for $f^*$; in particular, we may choose $f^*$ to be  right- continuous.\\
If $f^*$ is a version of $f\in L^1$ with minimal variation then Proposition \ref{Prop-2.2.3} applies to $f^*$ as well, with $\tau$ replaced by $\tau^{n_0}$. Thus, for all sufficiently large $n$,
\begin{equation*}
     \textrm{var}_I P_{\tau^{nn_0}} f^* \leq \frac{1}{2} \cdot \textrm{var}_I f^* +  B_n \cdot \|f^*\|_1,
\end{equation*}
and $ \textrm{var}_I  f^* = \textrm{var}_I f$. Since $P_{\tau^{nn_0}} f^*$ is a version of  $P_{\tau^{nn_0}} f$, we have
\begin{equation}\label{2.2.8}
    \begin{aligned}
        \textrm{var}_I P_{\tau^{nn_0}} f\leq \textrm{var}_I P_{\tau^{nn_0}} f^*\leq \frac{1}{2} \cdot \textrm{var}_I f^* +  B_{nn_0} \cdot \|f^*\|_1,\\
        = \frac{1}{2} \textrm{var}_I f+ B_{nn_0} \cdot \|f\|_1.
    \end{aligned}
\end{equation}
From equation (\ref{2.2.7}) and (\ref{2.2.8}) we get
\begin{equation*}
\begin{aligned}
     \|P_{\tau^{nn_0}} f\|_{BV} \leq \|f\|_1+ \frac{1}{2} \cdot \textrm{var}_I f+ B_{nn_0} \cdot \|f\|_1, \\
     = \frac{1}{2} \cdot \|f\|_{BV} +\left(\frac{1}{2}+B_{nn_0}\right) \|f\|_1 
     \end{aligned}
\end{equation*}
Choosing $C=\frac{1}{2}+B_{nn_0}$, we get the desired result.
\end{proof}
\begin{proposition}\label{prop 3.5}
\begin{enumerate}
    \item For every \( c > 0 \), the set \( E = \{ f \in L^1 : \| f \|_{\mathrm{BV}} \leq c \} \) is compact in \( L^1 \).
    \item \( (BV_I, \| \cdot \|_{\mathrm{BV}}) \) is a Banach space.
    \item \( BV_I \) is dense in \( L^1 \).
\end{enumerate}
For a detailed proof, see \cite{hofbauer1982}. 
\end{proposition}

\normalfont The properties of the operator $P_{\tau^n}$ and of the space $BV_I$ which we proved in Theorem \ref{Prop-2.2.3}  and in Propositions \ref{prop 2.2.5}, \ref{prop 3.5} allow us to use  Ionescu Tulcea and Marinescu theorem \cite{ionescu1950}, with $X=BV_I, Y=L^1$, and $P=P_{\tau^n}$ for the appropriate $n\in \mathbb{N}$.\\
Hofbauer and Keller \normalfont{\cites{hofbauer1982,keller1985} were the first to use this theorem for proving the quasi-compactness of $P_\tau$ and the existence of ACIM for $\tau$. Before we state our main result in this regard, we prove exactness of $\tau$ with ACIM.}\\
The proof is based on Theorem $3$ and  Theorem $4$ of \cite{lasota1982}.
\begin{definition}[Lower Function]
Let \( (I, \mathcal{B}, m) \) be a \(\sigma\)-finite measure space, and let \( \tau: I \to I \) be a nonsingular, doubly measurable transformation. A non-negative function \( h \in L^1\), \( h \geq 0\), \(\norm{h}_{L^1}>0\), is called a \emph{lower function} for \(P\) if
\begin{equation}\label{eq 15}
    \lim_{n \to \infty} \left\| (h - P^n f)^+ \right\|_{L^1} = 0,
\end{equation}
where \( (h - P^n f)^+ := \max\{h - P^n f, 0\} \) denotes the positive part of the function, and \( f \in \mathcal{D}(I, \mathcal{B}, m) \). It is enough if (\ref{eq 15}) holds for a dense subset of \( \mathcal{D}(I, \mathcal{B}, m) \). Here, \( \mathcal{D}(I, \mathcal{B}, m) \) denotes the set of all probability density functions on the measure space \( (I, \mathcal{B}, m) \).
\end{definition}
\begin{theorem}\label{thm 2.2.8}
    If \( \tau: I \to I \) is measurable and non-singular, and the operator \( P_{\tau} \) admits a lower function, then there exists a unique absolutely continuous invariant measure (ACIM) \( \mu \), and the system \( (\tau, \mu) \) is exact.
\end{theorem}
\begin{theorem}\label{thm 2.2.9} Let $\tau \in \mathcal{T}$. Then there exists the unique normalized absolutely continuous $\tau$ invariant measure $\mu$. The dynamical system $(I,\mathcal{B}, \mu, \tau)$ is exact and the density $\displaystyle h=\frac{d\mu}{dm}$ is bounded and decreasing.
\end{theorem}
\begin{proof}  The map $\tau$ satisfies conditions (2.1), (2.2) and (2.3) of Definition \ref{def 2.1.1}. We have proved in Proposition \ref{Prop 2.1.8} that $S$ is dense in $I$. Let $\chi_\Delta$ be the characteristic function of  an interval $\Delta=[d_0,d_1]$ whose end points belong to the set $S$. Now, we will prove that for all sufficiently large $n$, $P_\tau^n \chi_\Delta$ is a non-increasing function. We have proved that any iteration of $\tau$ satisfies the properties (2.1), (2.2) and (2.3). In particular, we have proved that $\tau^{(n)}_i$ is convex on $I_i^{(n)}=\left(a_i^{(n)},b_i^{(n)}\right)$, an element of the partition $\mathcal{P}^{(n)}$,  corresponding to $\tau^n$, and $\tau^n\left(a_i^{(n)}\right)=0$. This implies that
 $ g_n \cdot \left(\tau^{(n)}_i \right)^{-1}\cdot \chi_{\tau^{(n)}_i \left(a^{(n)}_i,b^{(n)}_i\right)}$
is a non-increasing function on $I$, since $g_n$ is non-increasing as the reciprocal of the derivative of a convex function.  We can see that,\\
\begin{equation*}
\left\{ a_1^{(n)},b_1^{(n)},a_2^{(n)},b_2^{(n)},\dots, a_i^{(n)},b_i^{(n)},\dots\right\}=\tau^{-n+1}\left\{a_1,b_1,a_2,b_2,\dots,a_i,b_i,\dots\right\}.
\end{equation*}
By the definition of the points $a_i^{(n)},b_i^{(n)}$, they are the preimages of the original partition points. This shows that $\mathcal{P}^{(n+1)}$ is a refinement of $\mathcal{P}^{(n)}$.
Since $d_1,d_2\in S$
there is an integer $n_0$ such that $d_1, d_2$ belong to the partition
$\left\{ a_1^{(n)},b_1^{(n)},a_2^{(n)},b_2^{(n)},\dots,a_i^{(n)},b_i^{(n)},\dots\right\}$ for $n\geq n_0$. The Frobenius-Perron operator for $\tau^n$ is
\begin{equation*}
    \begin{aligned}
        P_\tau^n f(x)=\sum_{y\in \tau^{-n}\{x\}} f(y) g_n(y).
    \end{aligned}
\end{equation*}
In particular, for $f=\chi_\Delta$ and $n\geq n_0$ we have,
\begin{equation*}
     P_\tau^n \chi_\Delta (x)=\sum_{y\in \tau^{-n}\{x\}} g_n(y) \cdot \chi_{\tau^n\bigl(I_i^{(n)}\bigl)}(y).
\end{equation*}
  Since $\tau^n\left(I_i^{(n)}\right)$ is of the form  $\left(0,\tau^n\left(b_i^{(n)}\right)\right)$, $P_\tau^n\chi_\Delta$ is non-increasing as a sum of non-increasing functions.
Now, let $D_0$ be the subset of $D(I, \mathcal{B}, m)$ consisting of all densities of the form
\begin{equation*}
    f(x)=\sum_{i=1}^N c_i \chi_{\Delta_i}(x), N\in \mathbb{N}, c_i\geq 0,
\end{equation*}
where the endpoints of the intervals $\Delta_i$ belong to $S$. Since S is dense in $I$, the set $D_0$ is dense in $D(I, \mathcal{B}, m)$. Now, we construct a lower function for $P_\tau$. Let $f\in D_0$. There exists $n_0=n_0(f)$ such that $P_\tau^n f$ is non-increasing for $n \geq n_0$. By part (1) of Theorem \ref{Thm 2.1.6} for any $\tau\in \mathcal{T}, P_\tau$ preserves the cone of non-increasing functions, and by (2), we have $P_\tau^n f(x)\leq 1/x$ for $n\geq n_0$.
Now, using this estimate and Theorem \ref{Thm 2.1.6} we get,
\begin{equation*}
\begin{aligned}
    P_\tau^{n+1}f(0)=P_\tau(P_\tau^n f(0))\leq  \alpha \cdot P_\tau^n f(0) + D,
\end{aligned}
\end{equation*}
where  $\alpha<1$ and $D$ are defined as in Theorem \ref{Thm 2.1.6} for both cases. Using an induction argument, we get
\begin{equation*}
     P_\tau^{n+n_0}f(0) \leq \alpha^n\cdot  P_\tau^{n_0} f(0)+\frac{D}{1-\alpha}.
\end{equation*}
Let $K=1+\frac{D}{1-\alpha}$. For sufficiently large $n$, say $n\geq n_1$, we have $P_\tau^n f(0) \leq K$.\\
Define $h=\frac{1}{2}\chi_{[0,1/(2K)]}$. We will prove that 
\begin{equation*}
    P_\tau^n f(x) \geq h(x) \text{ for $n\geq n_1$}.
\end{equation*}
Indeed, suppose there existed $x_0\in [0,1/(2K)]$ such that $ P_\tau^n f(x_0)< h(x_0)=\frac{1}{2}$. Then
\begin{equation*}
\begin{aligned}
    1=\int_{0}^{x_0} P_\tau^n f dx +\int_{x_0}^{1} P_\tau^n f dx\leq x_0 P_\tau^n f(0)+(1-x_0) P_\tau^n f(x_0)
   < \frac{1}{2K}\cdot K +\frac{1}{2} \hspace{.2 cm}=1,
    \end{aligned}
\end{equation*}
which is not possible. Hence $ P_\tau^n f(x) \geq h(x) \text{ for $n\geq n_1$}.$
\end{proof}
\normalfont{Theorem \ref{thm 2.2.9} implies that the only eigenvalue of $P_\tau$ of modulus $1$ is $1$, and that its eigenspace is one dimensional.   With Ionescu–Tulcea and Marinescu Theorem this gives the following :}
\theorem \label{Theorem 2.2.7} Let $\tau \in \mathcal{T}$. Then the Frobenius-Perron operator $P_\tau$ is quasi-compact on the space $BV_I$. More precisely, we have
\begin{enumerate}
\item $P_\tau: BV_I \rightarrow BV_I$ has $1$ as the only eigenvalue of modulus $1$;
\item  $E_1=\{f\in L^1 \mid P_\tau f=f\} \subseteq BV_I$, and $E_1$ is one-dimensional;
\item $\displaystyle P_\tau= \Psi +Q$, where $\Psi$ represents the projection onto the  eigenspace  $E_1$, $\|\Psi\|_1= 1$, and $Q$ is a linear operator on $L^1$, $\displaystyle \sup_{n\in \mathbf{N}} \|Q^n\|_1 < \infty$ and $Q \cdot \Psi = 0$;
\item $Q(BV_I)\subset BV_I$ and, considered as a linear operator on $(BV_I,\|.\|_{BV})$, $Q$ satisfies $\|Q^n\|_{BV} \leq H \cdot q^n$  $(n\geq 1)$ for some constants $H>0$ and $0<q<1$.
\end{enumerate}
\begin{proof} Since $P_\tau$ is a contraction on $L^1$, its powers are uniformly bounded. That is, we have $\|P_\tau^n\|_{1} \leq 1$ for all $n \in \mathbb{N}$.

From Proposition~\ref{prop 2.2.5}, we know that the operator satisfies a Lasota–Yorke inequality of the form $$\|P_\tau^n f\|_{BV} \leq \frac{1}{2} \|f\|_{BV} + C \|f\|_{1},
$$
for some constant $C > 0$ and all sufficiently large $n$. In addition, by Proposition~\ref{prop 3.5}, the embedding $
BV_I \hookrightarrow L^1
$ is compact, hence the operator $P_\tau: (BV_I, \|\cdot\|_{BV}) \to (L^1, \|\cdot\|_1)$ is compact. These properties verify the conditions required to apply the Ionescu–Tulcea and Marinescu Theorem.

Theorem~\ref{thm 2.2.9} further shows that the only eigenvalue of $P_\tau$ of modulus one is $\lambda = 1$, and that the associated eigenspace is one-dimensional and lies within $BV_I$. This establishes (1) and (2).
Finally, (3) and (4) follows from the spectral decomposition guaranteed by the 
Ionescu–Tulcea and Marinescu theorem, which states that the action of $P_\tau$ can be decomposed into a projection onto the finite-dimensional eigenspace and a remainder term whose powers decay exponentially in the operator norm on $BV_I$.
\end{proof}
\normalfont{Quasi-compactnes of $P_\tau$ implies several important ergodic properties of the system $(\tau,\mu)$, such as exponential decay of correlation, a Central Limit Theorem, and many other probabilistic consequences; see \cites{hofbauer1982,keller1985} for details.}
\begin{itemize}
\item \textbf{Weak Mixing:} Since $1$ is the only eigenvalue of $P_\tau$ with modulus $1$, and the corresponding eigenspace is one-dimensional, the system $(\tau, \mu)$ does not have any non-trivial periodic components. This implies that $(\tau, \mu)$ is weakly mixing and has several important statistical and ergodic properties, including:
\item \textbf{Exponential Decay of Correlations:} For functions of bounded variation, the correlation function decays exponentially fast. This means that for any two observables \( f, g \in BV_I \), there exist constants \( C > 0 \) and \( 0<q< 1 \) such that for all $n\in \mathbb{N}$,
    \[ 
    \left| \int_I f \cdot (g \circ \tau^n) \, d\mu - \int_I f \, d\mu \int_I g \, d\mu \right| \leq C \|f\|_{BV} \|g\|_{BV} q^n.
    \]
    \item \textbf{Central Limit Theorem:} The system satisfies the Central Limit Theorem, meaning the sum of observations (properly normalized) converges in distribution to a normal distribution. Specifically, for a function $f \in BV_I $ with $\int_I f d\mu = 0$, the sequence of partial sums \( S_n = \sum_{i=0}^{n-1} f \circ \tau^i \) satisfies
    \[ 
    \frac{S_n}{\sqrt{n}} \xrightarrow{d} \mathcal{N}\bigl(0, \sigma^2 \bigl),
    \]
    where \(\sigma^2\) is the variance given by
    \[ 
    \sigma^2 = \lim_{n \to \infty} \frac{1}{n} \int_I S_n^2 \, d\mu
    \]
and the convergence is in distribution.
    \item \textbf{Other Probabilistic Properties:}  $\mu$ is the equilibrium state for $\log g$ on $I$, i.e.,
    \begin{equation*}
        h(\mu)+\int \log g d\mu =\sup\left\{ h(\nu)+\int \log g d\nu : \text{ $\nu$ is a $\tau$-invariant probability on $I$ }\right\},
    \end{equation*}
        where $h(\nu)$ is the entropy of $(\tau,\nu)$.
\end{itemize}
\textbf{Note :}
     It is often believed that if a map \( \tau: [0,1] \to [0,1] \) is onto, then the corresponding invariant measure must be supported on the entire interval \( [0,1] \). However, this is not always true. This can be illustrated through an example of a piecewise linear, so also piecewise convex map:
$$\tau(x)=\begin{cases} 2x &\ \ \text{ if }\ x\in[0,1/4);\\
                          2x-1/2 & \ \ \text{ if }\ x\in[1/4,1/2);\\
                      2x-1 &\ \ \text{ if }\ x\in[1/2,1].
\end{cases}
$$
\begin{figure}[H]
    \centering
    \includegraphics[width=0.75\linewidth]{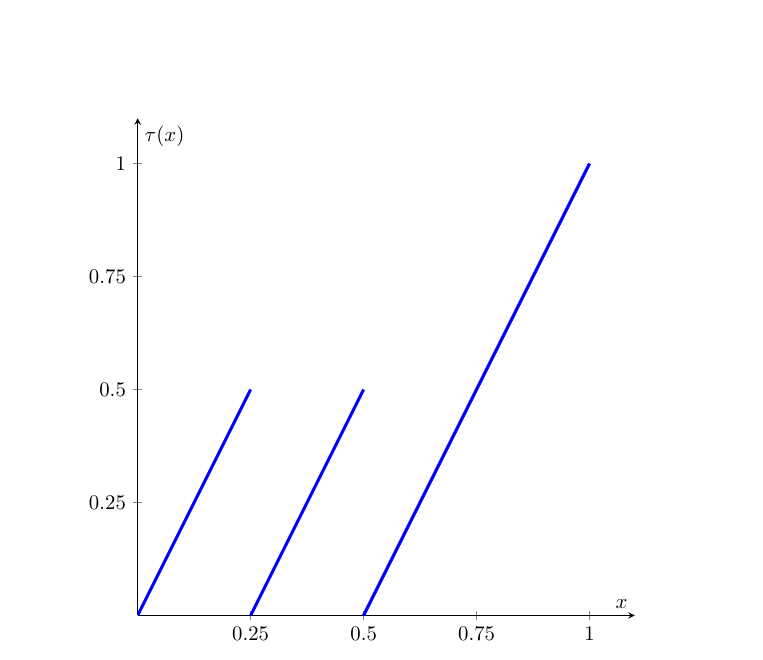}
    \caption{Plot of piecewise convex map $\tau(x)$.}
\end{figure}

While the range of the map covers the entire interval \( [0,1] \), the invariant density is
\begin{equation*}
    h(x)=\begin{cases} 2 &\ \ \text{ for }\ x\in[0,1/2);\\
                      0&\ \ \text{ for }\ x\in[1/2,1].
\end{cases}
\end{equation*}}

\end{document}